\documentclass[11pt,a4paper]{amsart} 
\usepackage[all]{xy}
\SelectTips{cm}{}
\usepackage{ifthen} 
\usepackage{amssymb}
\calclayout
\makeatletter 
\def\serieslogo@{} 
\def\@setcopyright{} 
\makeatother

\title{The Auslander-Reiten formula for complexes of modules}

\author{Henning Krause}
\address{Henning Krause\\ Institut f\"ur Mathematik\\
Universit\"at Paderborn\\ 33095 Paderborn\\ Germany}
\email{hkrause@math.uni-paderborn.de}
\author{Jue Le}
\address{Jue Le\\ Department of Mathematics, Beijing Normal University\\ 
100875 Beijing\\ P.R.China}
\email{lejue@vip.sina.com}

\thanks{Version from December 1, 2004. This work is supported by the
Doctor Program Foundation of the Ministry of Education of China (Grant
No. 20040027002)}

\newtheorem{lem}{Lemma}[section]
\newtheorem{prop}[lem]{Proposition} \newtheorem{cor}[lem]{Corollary}
\newtheorem{thm}[lem]{Theorem}

\theoremstyle{remark}

\theoremstyle{definition}
\newtheorem{exm}[lem]{Example}
\newtheorem{defn}[lem]{Definition}

\newtheorem{rem}[lem]{Remark}
\numberwithin{equation}{section}

\newcommand{\smatrix}[1]{\left[\begin{smallmatrix}#1\end{smallmatrix}\right]}

\renewcommand{\mod}{\operatorname{mod}\nolimits}

\newcommand{\proj}{\operatorname{proj}\nolimits}
\newcommand{\rad}{\operatorname{rad}\nolimits}

\newcommand{\id}{\operatorname{id}\nolimits}
\newcommand{\Tr}{\operatorname{Tr}\nolimits}

\newcommand{\Mod}{\operatorname{Mod}\nolimits}
\newcommand{\End}{\operatorname{End}\nolimits}
\newcommand{\Hom}{\operatorname{Hom}\nolimits}

\newcommand{\Ext}{\operatorname{Ext}\nolimits}

\newcommand{\umod}{\operatorname{\underline{mod}}\nolimits}
\newcommand{\uMod}{\operatorname{\underline{Mod}}\nolimits}
\newcommand{\oMod}{\operatorname{\overline{Mod}}\nolimits}

\newcommand{\uHom}{\operatorname{\underline{Hom}}\nolimits}
\newcommand{\oHom}{\operatorname{\overline{Hom}}\nolimits}

\newcommand{\Inj}{\operatorname{Inj}\nolimits}

\newcommand{\Proj}{\operatorname{Proj}\nolimits}

\newcommand{\Ab}{\mathrm{Ab}}
\newcommand{\ac}{\mathrm{ac}}

\newcommand{\op}{\mathrm{op}}
\newcommand{\inc}{\mathrm{inc}}

\newcommand{\can}{\mathrm{can}}
\newcommand{\Id}{\mathrm{Id}}

\newcommand{\comp}{\mathop{\raisebox{+.3ex}{\hbox{$\scriptstyle\circ$}}}}
\newcommand{\lto}{\longrightarrow}
\newcommand{\xto}{\xrightarrow}

\def\a{\alpha}
\def\b{\beta}

\def\d{\delta}
\def\g{\gamma}
\def\p{\phi}

\def\s{\sigma}

\def\Ga{\Gamma}
\def\La{\Lambda}

\def\A{{\mathcal A}}

\def\S{{\mathcal S}}

\def\T{{\mathcal T}}

\def\bfa{\mathbf a}

\def\bfi{\mathbf i}
\def\bfj{\mathbf j}
\def\bfp{\mathbf p}
\def\bft{\mathbf t}

\def\bfC{\mathbf C}
\def\bfD{\mathbf D}
\def\bfK{\mathbf K}

\begin{document}

\begin{abstract}
An Auslander-Reiten formula for complexes of modules is presented.
This formula contains as a special case the classical Auslander Reiten
formula. The Auslander-Reiten translate of a complex is described
explicitly, and various applications are discussed.
\end{abstract}
\maketitle

\section{Introduction}
The classical Auslander-Reiten formula for modules over a noetherian
algebra $\La$ says that
$$D\Ext_\La^1(M,N)\cong\oHom_\La(N,D\Tr M)$$ whenever $M$ is finitely
generated \cite{AR}. Here, $D$ denotes the duality over a fixed
commutative ground ring $k$. In this paper, we extend this to a
formula for complexes of modules. We take as morphisms chain maps up
to homotopy and obtain the formula for $D\Ext_\La^1(M,N)$ as a special
case by applying it to injective resolutions $\bfi M$ and $\bfi N$.

Various authors noticed the analogy between the Auslander-Reiten
formula and Serre duality for categories of sheaves; see for instance
\cite{Sc,GL}, and see \cite{DL} for a formulation of Serre duality in
terms of extension groups. Passing from abelian categories to their
derived categories, further analogies have been noticed, in particular
in connection with the existence of Auslander-Reiten triangles \cite{K1,RV}.

There is the following common setting for proving such duality
formulas. Let $\T$ be a $k$-linear triangulated category which is
compactly generated. Then one can apply Brown's representability
theorem and has for any compact object $X$ a representing object $\bft
X$ such that
$$D\Hom_\T(X,-)\cong\Hom_\T(-,\bft X).$$ In this paper, we take for $\T$
the category $\bfK(\Inj\La)$ of complexes of injective $\La$-modules
up to homotopy. Then we can prove that $\bft X=\bfp X\otimes_\La
D\La$, where $\bfp X$ denotes the projective resolution of $X$. 

There is a good reason to consider the category of complexes
$\bfK(\Inj\La)$. The injective resolutions of all finitely generated
modules generate the full subcategory of compact objects, which
therefore is equivalent to the bounded derived category
$\bfD^b(\mod\La)$ of the category $\mod\La$ of finitely generated
$\La$-modules.

Our identification of the translation $\bft $ has various interesting
consequences. For instance, we can describe Auslander-Reiten triangles
in $\bfK(\Inj\La)$, and we get a simple recipe for computing almost
split sequences in the category $\Mod\La$ of $\La$-modules which seems
to be new. 

There is another method for computing Auslander-Reiten triangles in
$\bfK(\Inj\La)$. This is based on the construction of an adjoint for
Happel's functor
$$\bfD^b(\mod\La)\lto\umod\hat\La$$
into the stable module category of
the repetitive algebra $\hat\La$ \cite{H}. To be precise, we extend
Happel's functor to a functor $\bfK(\Inj\La)\to\uMod\hat\La$ on
unbounded complexes, and this admits a right adjoint which preserves
Auslander-Reiten triangles.

The Auslander-Reiten theory for complexes has been initiated by
Happel. In \cite{H,H1}, he introduced Auslander-Reiten triangles and
characterized their existence in the derived category
$\bfD^b(\mod\La)$.  This pioneering work has been extended by various
authors. More recently, Bautista et al. established in \cite{BSZ} the
existence of almost split sequences in some categories of complexes of
fixed size.

The methods in this paper might be of interest not only for studying
module categories. In fact, they can be applied to other more general
settings. To illustrate this point, we include an Auslander-Reiten
formula for computing $\Ext^1_\A(-,-)$ in any locally noetherian
Grothendieck category $\A$.

\section{The homotopy category of injectives}

Let $k$ be a commutative noetherian ring which is complete and local.
Throughout this paper, we fix a noetherian $k$-algebra $\La$, that is,
a $k$-algebra which is finitely generated as a module over $k$.

We consider the category $\Mod\La$ of (right) $\La$-modules and 
the following full subcategories:
\begin{itemize}
\item[] $\mod\La$ = the finitely presented $\La$-modules,
\item[] $\Inj\La$ = the injective $\La$-modules,
\item[] $\Proj\La$ = the projective $\La$-modules,
\item[] $\proj\La$ = the finitely generated projective $\La$-modules.
\end{itemize}
Note that the assumptions on $\La$ imply that every finitely
generated $\La$-module decomposes essentially uniquely into a finite
coproduct of indecomposable modules with local endomorphism rings.

In addition, we fix an injective envelope $E=E(k/\mathfrak m)$, where
$\mathfrak m$ denotes the unique maximal ideal of $k$. We obtain a
duality
$$D=\Hom_k(-,E)\colon \Mod k\lto\Mod k$$ which induces a duality
between $\Mod\La$ and $\Mod\La^\op$.

Given any additive category $\A$, we denote by $\bfC(\A)$ the category
of cochain complexes in $\A$, and we write $\bfK(\A)$ for the category
of cochain complexes up to homotopy. If $\A$ is abelian, the derived
category is denoted by $\bfD(\A)$. We refer to \cite{V} for further
notation and terminology concerning categories of complexes.

We denote by  
$$\bfp\colon\bfD(\Mod\La)\lto\bfK(\Proj\La)$$ the left adjoint of
the composite
$$\bfK(\Proj\La)\xto{\inc}\bfK(\Mod\La)
\xto{\can}\bfD(\Mod\La)$$
which sends a complex $X$ to its {\em projective resolution} $\bfp X$.
Dually, we denote by  
$$\bfi\colon\bfD(\Mod\La)\lto\bfK(\Inj\La)$$ the right adjoint of
the composite
$$\bfK(\Inj\La)\xto{\inc}\bfK(\Mod\La) \xto{\can}\bfD(\Mod\La)$$ which
sends a complex $X$ to its {\em injective resolution} $\bfi X$. For
the existence of $\bfp$ and $\bfi$, see \cite{S,BN}.

We shall work in the category $\bfK(\Inj\La)$. This is a triangulated
category with arbitrary coproducts. We denote by $\bfK^c(\Inj\La)$ the
full subcategory which is formed by all compact objects. Recall that
an object $X$ is {\em compact} if every map $X\to\coprod_{i\in I} Y_i$
factors through $\coprod_{i\in J} Y_i$ for some finite $J\subseteq
I$. Let us collect from \cite{K2} the basic properties of $\bfK(\Inj\La)$.

\begin{lem}
An object in $\bfK(\Inj\La)$ is compact if and only if it is
isomorphic to a complex $X$ satisfying
\begin{enumerate}
\item $X^n=0$ for $n\ll 0$, 
\item $H^nX$ is finitely generated over $\La$ for
all $n$, and 
\item $H^nX=0$ for $n\gg 0$.
\end{enumerate}
\end{lem}

\begin{lem}
The canonical functor $\bfK(\Inj\La)\to\bfD(\Mod\La)$
induces an equivalence
$$\bfK^c(\Inj\La)\lto\bfD^b(\mod\La).$$ 
\end{lem}

\begin{lem}
The triangulated category
$\bfK(\Inj\La)$ is compactly generated, that is, it coincides with the
smallest full triangulated subcategory closed under all coproducts and
containing all compact objects.
\end{lem}

\section{The Auslander-Reiten formula}

In this section, the Auslander-Reiten formula for complexes is proved.
We begin with a number of simple lemmas. Given a pair of complexes
$X,Y$ of modules over $\La$ or $\La^\op$, we denote by $\Hom_\La(X,Y)$
and $X\otimes_\La Y$ the total Hom and the total tensor product
respectively, which are complexes of $k$-modules.

\begin{lem}\label{le:iso1}
Let $X,Y$ be complexes in $\bfC(\Mod\La)$. Then we
have in $\bfC(\Mod k)$ a natural map
\begin{equation}\label{eq:lem1}
Y\otimes_\La\Hom_\La(X,\La)\lto\Hom_\La(X,Y),
\end{equation}
which is an isomorphism if $X\in\bfC^-(\proj\La)$ and $Y\in\bfC^+(\Mod\La)$.
\end{lem}
\begin{proof}
  Given $\La$-modules $M$ and $N$, we have a map
  $$\s\colon N\otimes_\La\Hom_\La(M,\La)\lto\Hom_\La(M,N)$$
  which is
  defined by
$$\s(n\otimes\p)(m)=n\p(m).$$ This map is an isomorphism if $M$ is
finitely generated projective and extends to an isomorphism of
complexes provided $X$ and $Y$ are bounded in the appropriate direction.
\end{proof}

\begin{lem}\label{le:iso2}
Let $M,N$ be $\La$-modules and suppose that $M$ is finitely presented.
Then we have an isomorphism
\begin{equation}\label{eq:lem2}
M\otimes_\La\Hom_k(N,E)\cong\Hom_k(\Hom_\La(M,N),E).
\end{equation}
\end{lem}
\begin{proof}
  We have the isomorphism for $M=\La$ and therefore whenever $M$ has a
  presentation $\La^n\to\La^m\to M\to 0$, since
  $-\otimes_\La\Hom_k(N,E)$ and $\Hom_k(\Hom_\La(-,N),E)$ are both
  right exact.
\end{proof}

\begin{lem}\label{le:iso3}
Let $X,Y$ be complexes of $\La$-modules. Then we have an isomorphism
\begin{equation}\label{eq:lem3}
H^0\Hom_\La(X,Y)\cong\Hom_{\bfK(\Mod\La)}(X,Y).
\end{equation}
\end{lem}
\begin{proof} 
Taking cycles of $\Hom_\La(X,Y)$ in degree zero picks the chain maps
$X\to Y$, and boundaries in degree zero form the subgroup of
null-homotopic chain maps. Thus $H^0\Hom_\La(X,Y)$ equals the set of chain maps $X\to Y$ up to homotopy.
\end{proof}

Let us consider the following commutative diagram
$$\xymatrix{ &\bfD^-(\Mod\La)&\bfK^-(\Proj\La)\ar[l]_\sim\\
\bfK^c(\Inj\La)\ar[r]^\sim\ar@{>->}[d]&\bfD^b(\mod\La)\ar@{>->}[u]\ar@{>->}[d]&
\bfK^{-,b}(\proj\La)\ar[l]_\sim\ar@{>->}[u]\\
\bfK^+(\Inj\La)\ar[r]^\sim&\bfD^+(\Mod\La) }$$ in which all horizontal
functors are obtained by restricting the localization functor
$\bfK(\Mod\La)\to\bfD(\Mod\La)$ to appropriate subcategories.  We
denote by
$$\pi\colon\bfK^c(\Inj\La)\lto\bfK^{-,b}(\proj\La)$$ the composite of
the equivalence $\bfK^c(\Inj\La)\to\bfD^b(\mod\La)$ with a
quasi-inverse of the equivalence
$\bfK^{-,b}(\proj\La)\to\bfD^b(\mod\La)$. Note that $\pi X\cong \bfp X$.

\begin{thm}
Let $X$ and $Y$ be complexes of injective $\La$-modules.  Suppose that
$X^n=0$ for $n\ll 0$, that $H^nX$ is finitely generated over $\La$
for all $n$, and that $H^nX=0$ for $n\gg 0$.  Then we
have an isomorphism
\begin{equation}\label{eq:AR}
D\Hom_{\bfK(\Inj\La)}(X,Y)\cong
\Hom_{\bfK(\Inj\La)}(Y,(\pi X)\otimes_\La D\La)
\end{equation}
which is natural in $X$ and $Y$.
\end{thm}
\begin{proof}
We use the fact that $\bfK(\Inj\La)$ is compactly generated. Therefore
it is sufficient to verify the isomorphism for every compact object
$Y$. This follows from the subsequent Lemma~\ref{le:inj}.  Thus we
suppose that $Y$ is a compact object in $\bfK(\Inj\La)$. Note that
this implies $Y^n=0$ for $n\ll 0$, and in particular $Y\cong\bfi
Y$. We obtain the following sequence of isomorphisms, where
short arguments are added on the right hand side.
\begin{alignat*}{2}
D\Hom_{\bfK(\Inj\La)}(X,Y)
&\cong\Hom_k(\Hom_{\bfK(\Inj\La)}(X,\bfi Y),E)&&\qquad\text{$Y$ compact}\\
&\cong\Hom_k(\Hom_{\bfD(\Mod\La)}(X,Y),E)&&\qquad\text{adjunction}\\
&\cong\Hom_k(\Hom_{\bfK(\Mod\La)}(\pi X,Y),E)&&\qquad\text{adjunction}\\
&\cong\Hom_k(H^0\Hom_{\La}(\pi X,Y),E)&&\qquad\text{from (\ref{eq:lem3})}\\
&\cong H^0\Hom_k(\Hom_{\La}(\pi X,Y),E)&&\qquad\text{$E$ injective}\\
&\cong H^0\Hom_k(Y\otimes_\La\Hom_\La(\pi X,\La),E)
&&\qquad\text{from (\ref{eq:lem1})}\\
&\cong H^0\Hom_\La(Y,\Hom_k(\Hom_\La(\pi X,\La),E))&&\qquad\text{adjunction}\\
&\cong H^0\Hom_\La(Y,(\pi X)\otimes_\La\Hom_k(\La,E))
&&\qquad\text{from (\ref{eq:lem2})}\\
&\cong \Hom_{\bfK(\Inj\La)}(Y,(\pi X)\otimes_\La D\La)
&&\qquad\text{from (\ref{eq:lem3})}
\end{alignat*}
This isomorphism completes the proof.
\end{proof}

\begin{lem}\label{le:inj}
Let $X,X'$ be objects in a $k$-linear compactly generated triangulated
category. Suppose that $X$ is compact. If there is a natural
isomorphism
$$D\Hom_\T(X,Y)\cong\Hom_\T(Y,X')$$ for all compact $Y\in\T$, then
$D\Hom_\T(X,-)\cong\Hom_\T(-,X')$.
\end{lem}
\begin{proof}
We shall use Theorem~1.8 in \cite{K0}, which
states the following equivalent conditions for an object $Y$ in $\T$.
\begin{enumerate}
\item The object $Y$ is pure-injective.
\item The object $H_Y=\Hom_\T(-,Y)|_{\T^c}$ is injective in the
category $\Mod\T^c$ of contravariant additive functors $\T^c\to\Ab$.
\item The map $\Hom_\T(Y',Y)\to\Hom_{\T^c}(H_{Y'},H_Y)$ sending
$\p$ to $H_\p$ is bijective for all $Y'$ in $\T$.
%\item If $H_\p=0$ for some map $\p\colon Y'\to Y$ in $\T$, then $\p=0$.
\end{enumerate}
Here, $\T^c$ denotes the full subcategory of compact objects in $\T$.

We apply Brown's representability theorem (see \cite[5.2]{K} or
\cite[Theorem~3.1]{N1}) and obtain  an object $X''$ such that
$$D\Hom_\T(X,-)\cong\Hom_\T(-,X''),$$
since $X$ is compact. Condition (2) implies that both
objects $X'$ and $X''$ are pure-injective, since
$\Hom_{\T^c}(X,-)$ is a projective object in the category of
covariant additive functors $\T^c\to\Ab$, by Yoneda's lemma. We
have an isomorphism
$$H_{X'}=\Hom_\T(-,X')|_{\T^c}\cong
D\Hom_\T(X,-)|_{\T^c}\cong\Hom_\T(-,X'')|_{\T^c}=H_{X''},$$ and (3)
implies that this isomorphism is induced by an isomorphism $X'\to X''$
in $\T$. We conclude that
$$D\Hom_\T(X,-)\cong\Hom_\T(-,X').$$
\end{proof}

\section{The Auslander-Reiten translation}
In this section, we investigate the properties of the Auslander-Reiten
translation for complexes of $\La$-modules. The Auslander-Reiten
translation $D\Tr$ for modules is obtained from the translation for
complexes.  In particular, we deduce the classical Auslander-Reiten
formula.

\begin{prop} 
The functor
$$\bft \colon\bfK(\Inj\La)\xto{\can}\bfD(\Mod\La)
\xto{\bfp}\bfK(\Proj\La)\xto{-\otimes_\La D\La}\bfK(\Inj\La)$$ has the
following properties.
\begin{enumerate}
\item $\bft $ is exact and preserves all coproducts.
\item For compact objects
$X,Y$ in $\bfK(\Inj\La)$, the natural map $$\Hom_{\bfK(\Inj\La)}(X,Y)\lto
\Hom_{\bfK(\Inj\La)}(\bft X,\bft Y)$$ is bijective. 
\item For $X,Y$ in $\bfK(\Inj\La)$ with $X$ compact, there is a natural
isomorphism 
$$D\Hom_{\bfK(\Inj\La)}(X,Y)\cong \Hom_{\bfK(\Inj\La)}(Y,\bft X).$$
\item $\bft$ admits a right adjoint which is $\bfi\Hom_\La(D\La,-)$.
\end{enumerate}
\end{prop}
\begin{proof} 
(1) is clear and (3) follows from (\ref{eq:AR}). 
Now observe that for each pair $X,Y$ of compact objects,
the $k$-module $\Hom_{\bfK(\Inj\La)}(X,Y)$ is finitely generated.
Therefore (2) follows from (3), since we have the isomorphism
$$\begin{aligned}
\Hom_{\bfK(\Inj\La)}(X,Y)&\cong D^2\Hom_{\bfK(\Inj\La)}(X,Y)\\ &\cong
D\Hom_{\bfK(\Inj\La)}(Y,\bft X)\\ &\cong \Hom_{\bfK(\Inj\La)}(\bft X,\bft Y).
\end{aligned}$$
To prove (4), let $X,Y$ be objects in $\bfK(\Inj\La)$.
Then we have
$$\begin{aligned}
\Hom_{\bfK(\Mod\La)}(\bfp X\otimes_\La D\La,Y)&\cong
\Hom_{\bfK(\Mod\La)}(\bfp X,\Hom_\La (D\La,Y))\\
&\cong \Hom_{\bfD(\Mod\La)}(X,\Hom_\La (D\La,Y))\\
&\cong \Hom_{\bfK(\Mod\La)}(X,\bfi \Hom_\La (D\La,Y)).
\end{aligned}$$ 
Thus $\bft$ and $\bfi\Hom_\La(D\La,-)$ form an adjoint pair.
\end{proof}

Let us continue with some definitions. We denote by $\uMod\La$ the
{\em stable module category} modulo projectives which is obtained by
forming for each pair of $\La$-modules $M$ and $N$ the quotient
$$\uHom_\La(M,N)=\Hom_\La(M,N)/\{M\to P\to N\mid P
\textrm{ projective}\}.$$ Analogously, the stable module category
$\oMod\La$ modulo injectives is defined. 

Recall that a $\La$-module $M$ is {\em finitely presented} if it
admits a projective presentation $$P_1\xto{} P_0\to M\to 0$$ such that
$P_0$ and $P_1$ are finitely generated.  The {\em transpose} $\Tr M$
relative to this presentation is the $\La^\op$-module which is defined
by the exactness of the induced sequence
$$\Hom_\La(P_0,\La)\xto{}\Hom_\La(P_1,\La)\to\Tr M\to 0.$$ Note that
the presentation of $M$ is minimal if and only if the corresponding
presentation of $\Tr M$ is minimal. The construction of the
transpose is natural up to maps factoring through a projective and
induces a duality $\umod\La\to\umod\La^\op$.

\begin{prop}
The functor
$$\bfa\colon \Mod\La\xto{\inc}\bfD(\Mod\La)\xto{\bfi}\bfK(\Inj\La)\xto{\bft}
\bfK(\Inj\La)\xto{Z^{-1}}\oMod\La$$  has the
following properties.
\begin{enumerate}
\item $\bfa M\cong D\Tr M$ for every finitely presented $\La$-module $M$.
\item $\bfa $ preserves all coproducts.
\item $\bfa $ annihilates all projective $\La$-modules and induces a functor
$\uMod\La\to\oMod\La$.
\item Each exact sequence $0\to L\to M\to N\to 0$ of $\La$-modules
induces a sequence
$$0\to \bfa L\to \bfa M\to \bfa N\to L\otimes_\La D\La\to M\otimes_\La
D\La\to N\otimes_\La D\La\to 0$$ of $\La$-modules which is exact.
\end{enumerate}
\end{prop}
\begin{proof}
(1) The functor $\bft $ sends an injective resolution $\bfi M$ of $M$ to $\bfp
M\otimes_\La D\La$.  Using (\ref{eq:lem2}), we have
$$\bfp M\otimes_\La D\La\cong D\Hom_\La(\bfp M,\La).$$ This implies
\begin{equation}\label{eq:DTr}
Z^{-1}(\bfp M\otimes_\La D\La)\cong D\Tr M.
\end{equation}

(2) First observe that $\coprod_i(\bfi M_i)\cong\bfi (\coprod_i M_i)$
for every family of $\La$-modules $M_i$.  Clearly, $\bft $ and $Z^{-1}$
preserve coproducts. Thus $\bfa $ preserves coproducts.

(3) We have $\bft (\bfi \La)=D\La$ and therefore $\bfa \La=0$ in
$\oMod\La$. Thus $\bfa $ annihilates all projectives since it preserves
coproducts.

(4) An exact sequence $0\to L\to M\to N\to 0$ induces an exact
triangle $\bfp L\to \bfp M\to \bfp N\to (\bfp L)[1]$. This triangle
can be represented by a sequence $0\to \bfp L\to \bfp M\to \bfp N\to 0$ of
complexes which is split exact in each degree. Now apply $-\otimes_\La
D\La$ and use the snake lemma.
\end{proof}

We are now in the position to deduce the classical Aulander-Reiten
formula for modules \cite{AR} from the formula for complexes.

\begin{cor}[Auslander/Reiten]
Let $M$ and $N$ be $\La$-modules and suppose that $M$ is finitely
presented. Then we have an isomorphism
\begin{equation}\label{eq:carf}
D\Ext_\La^1(M,N)\cong\oHom_\La(N,D\Tr M).
\end{equation}
\end{cor}
\begin{proof}
Let $\bfi M$ and $\bfi N$ be injective resolutions of $M$ and $N$,
respectively. We apply the Auslander-Reiten formula (\ref{eq:AR}) and
the formula (\ref{eq:DTr}) for the Auslander-Reiten translate. Thus we
have
$$ D\Ext^1_\La(M,N)\cong D\Hom_{\bfK(\Inj\La)}(\bfi M,(\bfi
N)[1])\cong \Hom_{\bfK(\Inj\La)}(\bfi N,(\bfp M\otimes_\La
D\La)[-1]),$$ and the map
$$\Hom_{\bfK(\Inj\La)}(\bfi N,(\bfp M\otimes_\La
D\La)[-1])\lto\oHom_\La(N,D\Tr M),\quad \p\mapsto Z^0\p,$$ is clearly
surjective. The composite is also an injective map, since a map $N\to
D\Tr M$ factoring through an injective module $N'$ comes from an
element in $D\Ext^1_\La(M,N')$ which vanishes.
\end{proof}

\section{A general Auslander-Reiten formula for $\Ext_\A^1(-,-)$}

In this section, we extend the classical Auslander-Reiten formula for
modules to a formula for a more general class of abelian
categories. Thus we fix a locally noetherian Grothendieck category
$\A$, that is, $\A$ is an abelian Grothendieck category having a set
of generators which are noetherian objects in $\A$.

\begin{thm}
Let $M$ and $N$ be objects in $\A$.  Suppose that $M$ is noetherian and let
$\Ga=\End_\A(M)$. Given an injective $\Ga$-module $I$, there is an
object $\bfa_IM$ in $\A$ and an isomorphism
$$\Hom_\Ga(\Ext^1_\A(M,N),I)\cong\oHom_\A(N,\bfa_IM)$$
which is natural in $I$ and $N$.
\end{thm}
\begin{proof}
The category $\bfK(\Inj\A)$ is compactly generated and an injective
resolution $\bfi M$ of $M$ is a compact object; see
\cite[Proposition~2.3]{K2}. The functor
$$\Hom_\Ga(\Hom_{\bfK(\Inj\A)}(\bfi M,-),I)$$ is cohomological and sends
coproducts in $\bfK(\Inj\A)$ to products of abelian groups.  Using
Brown representability (see \cite[5.2]{K} or \cite[Theorem~3.1]{N1}),
we have a representing object $\bft_I(\bfi M)$ in $\bfK(\Inj\A)$ such
that
$$\Hom_\Ga(\Hom_{\bfK(\Inj\A)}(\bfi
M,-),I)\cong\Hom_{\bfK(\Inj\A)}(-,\bft_I(\bfi M)).$$
Now put
$\bfa_IM=Z^0\bft_I(\bfi M)$ and adapt the proof of the classical
Auslander-Reiten formula (\ref{eq:carf}) from the previous section.
\end{proof}

\begin{rem} 
If $\A$ is a $k$-linear category, then we can take $I=D\Ga$ and obtain
$$D\Ext^1_\A(M,N)\cong\oHom_\A(N,\bfa_{D\Ga}M).$$
\end{rem}

\section{Auslander-Reiten triangles}

In \cite{AR}, Auslander and Reiten used the formula
$$D\Ext_\La^1(M,N)\cong\oHom_\La(N,D\Tr M)$$ to establish the existence of
almost split sequences. More precisely, for a finitely presented
indecomposable and non-projective $\La$-module $M$, there exists an
almost split sequence
$$0\to D\Tr M\to L\to M\to 0$$ in the category of $\La$-modules.

In this section, we produce Auslander-Reiten triangles in the category
$\bfK(\Inj\La)$, using the Auslander-Reiten formula for complexes. In
addition, we show that almost split sequences can be obtained from
Auslander-Reiten triangles. This yields a simple recipe for the
construction of an almost split sequence.

Let us recall the relevant definitions from Auslander-Reiten theory. A
map $\a\colon X\to Y$ is called {\em left almost split}, if $\a$ is not
a section and if every map $X\to Y'$ which is not a section factors
through $\a$. Dually, a map $\b\colon Y\to Z$ is {\em right
almost split}, if $\b$ is not a retraction and if every map $Y'\to Z$
which is not a retraction factors through $\b$.

\begin{defn}
\begin{enumerate} 
\item An exact sequence $0\to X\xto{\a} Y\xto{\b} Z\to
0$ in an abelian category is called {\em almost split sequence}, if
$\a$ is left almost split and $\b$ is right almost split.
\item An exact triangle $X\xto{\a} Y\xto{\b}
Z\xto{\g} X[1]$ in a triangulated category is called {\em
Auslander-Reiten triangle}, if $\a$ is left almost split and $\b$ is
right almost split.
\end{enumerate}
\end{defn}

Happel introduced Auslander-Reiten triangles and studied their
existence in $\bfD^b(\mod\La)$ \cite{H}.  There is a general existence
result for Auslander-Reiten triangles in compactly generated
triangulated categories; see \cite{K1}. This yields the following.

\begin{prop}\label{pr:existence} 
Let $Z$ be a compact object in $\bfK(\Inj\La)$ which is
indecomposable.  Then there exists an Auslander-Reiten triangle
\begin{equation}\label{eq:art}
(\bfp Z\otimes_\La D\La)[-1]\xto{\a} Y\xto{\b} Z\xto{\g} \bfp
Z\otimes_\La D\La.
\end{equation}
\end{prop}
\begin{proof}
First observe that $\Gamma=\End_{\bfK(\Inj\La)}(Z)$ is local because
it is a noetherian $k$-algebra.  Now apply Theorem~2.2 from \cite{K1}
and use the formula (\ref{eq:AR}). Note that $D=\Hom_k(-,E)$ is
isomorphic to the functor $\Hom_\Gamma(-,I)$ with
$I=E(\Gamma/\rad\Gamma)$ which is used in \cite{K1}.
\end{proof}
Let us mention that the exact triangle (\ref{eq:art}) is determined by
the map $$\g\colon Z\lto \bfp Z\otimes_\La D\La$$ which corresponds under
the isomorphism
$$D\Hom_{\bfK(\Inj\La)}(Z,Z)\cong\Hom_{\bfK(\Inj\La)}(Z,\bfp
Z\otimes_\La D\La)$$ to a non-zero map $\End_{\bfK(\Inj\La)}(Z)\to E$
annihilating the radical of $\End_{\bfK(\Inj\La)}(Z)$.

An Auslander-Reiten triangle ending in the injective resolution of a
finitely presented indecomposable non-projective module induces an
almost split sequence as follows.

\begin{thm}\label{th:ass}
Let $N$ be a finitely presented $\La$-module which is indecomposable
and non-projective. Then there exists an Auslander-Reiten triangle
$$(\bfp N\otimes_\La D\La)[-1]\xto{\a} Y\xto{\b} \bfi N\xto{\g} \bfp
N\otimes_\La D\La$$ in $\bfK(\Inj\La)$ which the functor $Z^0$ sends
to an almost split sequence
\begin{equation*}\label{eq:ass}
0\to D\Tr N\xto{Z^0\a} Z^0Y\xto{Z^0\b} N\to 0
\end{equation*} in the
category of $\La$-modules.
\end{thm}
\begin{proof}
  The Auslander-Reiten triangle for $\bfi N$ is obtained from the
  triangle (\ref{eq:art}) by taking $Z=\bfi N$. Let us assume that the
  projective resolution $\bfp N$ is minimal. Note that we have a
  sequence 
\begin{equation}\label{eq:comp}
0\to (\bfp N\otimes_\La D\La)[-1]\xto{\a} Y\xto{\b} \bfi N\to 0
\end{equation} 
of chain maps which is split exact in each degree. This sequence is
  obtained from the mapping cone construction for $\g\colon\bfi
  N\xto{} \bfp N\otimes_\La D\La$.

We know from (\ref{eq:DTr}) that
$$D\Tr N\cong Z^0(\bfp N\otimes_\La D\La)[-1].$$ This module is
indecomposable and has a local endomorphism ring. Here we use that $N$
is indecomposable and that the resolution $\bfp N$ is minimal.  The
functor $Z^0$ takes the sequence (\ref{eq:comp}) to an exact sequence
\begin{equation}\label{eq:asss}
0\to D\Tr N\xto{Z^0\a} Z^0Y\xto{Z^0\b} N.
\end{equation}
Now observe that the map $Z^0\b$ is right almost split. This is clear
because $\b$ is right almost split and $Z^0$ induces a bijection
$\Hom_{\bfK(\Inj\La)}(\bfi M,\bfi N)\to\Hom_\La(M,N)$ for all $M$.  In
particular, $Z^0\b$ is an epimorphism since $N$ is non-projective. We
conclude from the following Lemma~\ref{le:Aus} that the sequence
(\ref{eq:asss}) is almost split.
\end{proof}

\begin{lem}\label{le:Aus}
An exact sequence $0\to X\xto{\a} Y\xto{\b} Z\to 0$ in an abelian
category is almost split if and only if $\b$ is right almost split and
the endomorphism ring of $X$ is local.
\end{lem}
\begin{proof} See Proposition~II.4.4 in \cite{A}. 
\end{proof}

\begin{rem} 
There is an analogue of Theorem~\ref{th:ass} for a projective module
$N$. Then $D\Tr N=0$ and $Z^0\b$ is the right almost split map ending in $N$.
\end{rem}

From the mapping cone construction for complexes, we get an explicit
recipe for the construction of an almost split sequence. Note that the
computation of almost split sequences is a classical problem in
representation theory \cite{G}. In particular, the middle term is
considered to be mysterious.

\begin{cor}
Let $N$ be a finitely presented $\La$-module which is indecomposable
and non-projective. Denote by
$$P_1\xto{\d_1} P_0\to N\to 0\quad\textrm{and}\quad 0\to N\to
I^0\xto{\d^0}I^1$$ a minimal projective presentation and an injective
presentation of $N$ respectively.  Choose a non-zero $k$-linear map
$\End_\La(N)\to E$ annihilating the radical of $\End_\La(N)$, and
extend it to a $k$-linear map $\p\colon\Hom_\La(P_0,I^0)\to E$. Let
$\bar\p$ denote the image of $\p$ under the isomorphism
$$D\Hom_\La(P_0,I^0)\cong\Hom_\La(I^0,P_0\otimes_\La D\La).$$
Then we have a commutative diagram with exact rows and columns
$$\xymatrix{&0\ar[d]&0\ar[d]&0\ar[d]\\ 
0\ar[r]& L\ar[r]^-{\a}\ar[d]& M\ar[r]^-{\b}\ar[d]&N\ar[r]\ar[d]&0\\ 
0\ar[r]& P_1\otimes_\La D\La\ar[r]^-{\smatrix{1\\0}}\ar[d]^-{\d_1\otimes 1}
& (P_1\otimes_\La D\La)\amalg I^0\ar[r]^-{\smatrix{0&1}}
\ar[d]^-{\smatrix{\d_1\otimes 1&\bar\p\\0&\d^0}}&I^0\ar[r]\ar[d]^-{\d^0}&0\\
0\ar[r]& P_0\otimes_\La D\La\ar[r]^-{\smatrix{1\\0}}& 
(P_0\otimes_\La D\La)\amalg I^1\ar[r]^-{\smatrix{0&1}}&I^1\ar[r]&0}$$
such that 
%$$0\to D\Tr N\xto{\a} M\xto{\b} N\to 0$$ 
the upper row is an almost split sequence
in the category of $\La$-modules.
\end{cor}

\begin{exm}
Let $k$ be a field and $\La=k[x]/(x^2)$.  Let $\bfi S$ denote the
injective resolution of the unique simple $\La$-module $S=k[x]/(x)$.
The corresponding Auslander-Reiten triangle in $\bfK(\Inj\La)$ has the
form
$$\bfp S[-1]\xto{\a} Y\xto{\b}\bfi S\xto{\g}\bfp S,$$ where $\g$
denotes an arbitrary non-zero map.  Viewing $\La$ as a complex
concentrated in degree zero, the corresponding Auslander-Reiten
triangle has the form
$$\La[-1]\xto{\a} Y\xto{\b}\La\xto{\g}\La,$$ where
$\g$ denotes the map induced by multiplication with $x$.
\end{exm}

\section{An adjoint of Happel's functor}

Let $\La$ be an artin $k$-algebra, that is, we assume that $\La$ is
artinian as a module over $k$. We denote by $\hat\La$ its repetitive
algebra. In this section, we extend Happel's functor \cite{H}
$$\bfD^b(\mod\La)\lto\umod\hat\La$$ to a functor which is defined on
unbounded complexes, and we give a right adjoint.

The repetitive algebra is by definition the doubly infinite matrix algebra
without identity
\setlength{\unitlength}{12pt}
$$\hat{\La}=\smatrix{
  \makebox(1,1){$\scriptstyle{\ddots}$}&&&&\\
  \makebox(1,1){$\scriptstyle{\ddots}$}&\makebox(1,1){$\scriptstyle{\La}$}
  &&&\makebox(1,1){$0$}\\
  &\makebox(1,1){$\scriptstyle{D\La}$}&\makebox(1,1){$\scriptstyle{\La}$}\\
  &&\makebox(1,1){$\scriptstyle{D\La}$}&\makebox(1,1){$\scriptstyle{\La}$}\\
  \makebox(1,1){$0$}&&&\makebox(1,1){$\scriptstyle{\ddots}$}&\makebox(1,1){$\scriptstyle{\ddots}$}\\ \\
}$$
in which matrices have only finitely many non-zero entries and the
multiplication is induced from the canonical maps $\La \otimes_\La
D\La \to D\La$ , $D\La \otimes_\La \La \to D\La$ , and the zero map
$D\La \otimes_\La D\La \to 0$. Note that projective and injective
modules over $\hat\La$ coincide. We denote by $\bfK_\ac(\Inj\hat\La)$
the full subcategory of $\bfK(\Inj\hat\La)$ which is formed by all
acyclic complexes. The follwing description of the stable category
$\uMod\hat\La$ is well-known; see for instance \cite[Example~7.6]{K2}.

\begin{lem}
  The functor $Z^0\colon\bfK_\ac(\Inj\hat\La)\to\uMod\hat\La$ is an
  equivalence of triangulated categories.
\end{lem}

We consider the algebra homomorphism
$$\p\colon \hat\La\lto\La, \quad(x_{ij})\mapsto x_{00},$$ and we view
$\La$ as a bimodule $_\La\La_{\hat\La}$ via $\p$. Let us explain the
following diagram.
$$\xymatrix{ \mod\La\ar[rr]^{-\otimes_\La\La}\ar[d]^\inc&&
\mod\hat\La\ar[d]^\inc\ar[rr]^\can&&\umod\hat\La\ar@{=}[d]\\
\bfD^b(\mod\La)\ar[rr]^{-\otimes_\La\La}&&\bfD^b(\mod\hat\La)\ar[rr]&&
\umod\hat\La\\
\bfK^c(\Inj\La)\ar[rr]^{F^c}\ar[d]^\inc\ar[u]_-\wr&&\bfK^c(\Inj\hat\La)\ar[rr]^{G^c}\ar[d]^\inc\ar[u]_-\wr&&
\bfK_\ac^c(\Inj\hat\La)\ar[d]^\inc\ar[u]^-{Z^0}_-\wr\\
\bfK(\Inj\La)\ar@<1ex>[rr]^F\ar@<1ex>[d]^\inc&&\bfK(\Inj\hat\La)\ar@<1ex>[ll]^{\Hom_{\hat\La}(\La,-)}\ar@<1ex>[d]^\inc\ar@<1ex>[rr]^G&&
\bfK_\ac(\Inj\hat\La)\ar@<1ex>[ll]^-\inc\\
\bfK(\Mod\La)\ar@<1ex>[u]^{\bfj_\La}\ar@<1ex>[rr]^{-\otimes_\La\La}&&
\bfK(\Mod\hat\La)\ar@<1ex>[u]^{\bfj_{\hat\La}}\ar@<1ex>[ll]^{\Hom_{\hat\La}(\La,-)}}$$
The top squares show the construction of Happel's functor
$\bfD^b(\mod\La)\to\umod\hat\La$ for which we refer to \cite[2.5]{H1}.

The bimodule $_\La\La_{\hat\La}$ induces an adjoint pair of functors
between $\bfK(\Mod\La)$ and $\bfK(\Mod\hat\La)$. Note that
$\Hom_{\hat\La}(\La,-)$ takes injective $\hat\La$-modules to injective
$\La$-modules. Thus we get an induced functor
$\bfK(\Inj\hat\La)\to\bfK(\Inj\La)$. This functor preserves products
and has therefore a left adjoint $F$, by Brown's representability
theorem \cite[Theorem~8.6.1]{N2}. A left adjoint preserves compactness
if the right adjoint preserves coproducts; see \cite[Theorem~5.1]{N1}.
Clearly, $\Hom_{\hat\La}(\La,-)$ preserves coproducts since $\La$ is
finitely generated over $\hat\La$.  Thus $F$ induces a functor $F^c$.

The inclusion $\bfK(\Inj\La)\to\bfK(\Mod\La)$ preserves products and
has therefore a left adjoint $\bfj_\La$, by Brown's representability
theorem \cite[Theorem~8.6.1]{N2}. Note that $\bfj_\La M=\bfi M$ is an
injective resolution for every $\La$-module $M$. We have the same for
${\hat\La}$, of course. Thus we have
$$F\comp\bfj_\La=\bfj_{\hat\La}\comp (-\otimes_\La\La).$$ It follows
that $F$ takes the injective resolution of a $\La$-module $M$ to the
injective resolution of the $\hat\La$-module $M\otimes_\La\La$. This
shows that $F^c$ coincides with $-\otimes_\La\La$ when one passes to
the derived category $\bfD^b(\mod\La)$ via the canonical equivalence
$\bfK^c(\Inj\La)\to\bfD^b(\mod\La)$.

The inclusion $\bfK_\ac(\Inj\hat\La)\to\bfK(\Inj\hat\La)$ has a left
adjoint $G$; see \cite[Theorem~4.2]{K2}. This left adjoint admits
an explicit description. For instance, it takes the injective
resolution $\bfi M$ of a $\hat\La$-module $M$ to the mapping cone of
the canonical map $\bfp M\to\bfi M$, which is a complete resolution of
$M$. The functor $G$ preserves compactness and induces therefore a
functor $G^c$, because its right adjoint preserves coproducts
\cite[Theorem~5.1]{N1}.

The following result summarizes our construction.

\begin{thm}
The composite 
$$\uMod\hat\La\xto{\sim}\bfK_\ac(\Inj\hat\La)\xto{\Hom_{\hat\La}(\La,-)}
\bfK(\Inj\La)$$ has a fully faithful left adjoint
$$\bfK(\Inj\La)\xto{G\comp
F}\bfK_\ac(\Inj\hat\La)\xto{\sim}\uMod\hat\La$$ which extends Happel's
functor
$$\bfD^b(\mod\La)\xto{-\otimes_\La\La}\bfD^b(\mod\hat\La)\to\umod\hat\La.$$
\end{thm}

\section{The computation of Auslander-Reiten triangles}

In this section, we explain a method for computing Auslander-Reiten
triangles in $\bfK(\Inj\La)$. It is shown that the adjoint of Happel's
functor reduces the computation to the problem of computing almost
split sequences in $\mod\hat\La$. This is based on the following result.

\begin{prop}\label{pr:adj}
Let $F\colon\S\to\T$ be a fully faithful exact functor between
triangulated categories which admits a right adjoint $G\colon\T\to\S$.
Suppose $$
X_\S\xto{\a_\S} Y_\S\xto{\b_\S}Z_\S\xto{\g_\S} X_\S[1]\quad\textrm{and}\quad
X_\T\xto{\a_\T} Y_\T\xto{\b_\T}Z_\T\xto{\g_\T} X_\T[1]
$$ are Auslander-Reiten triangles in $\S$ and $\T$ respectively, where
$Z_\T=FZ_\S$. Then $$GX_\T\xto{G\a_\T}
GY_\T\xto{G\b_\T}GZ_\T\xto{G\g_\T} GX_\T[1]$$ is the coproduct of
$X_\S\xto{\a_\S} Y_\S\xto{\b_\S}Z_\S\xto{\g_\S} X_\S[1]$ and a triangle
$W\xto{\id}W\to 0\to W[1]$.
\end{prop}
\begin{proof}
We have a natural isomorphism $\Id_\S\cong G\comp F$ which we view as
an identification. In particular, $G$ induces a bijection
\begin{equation}\label{eq:G}
\Hom_\T(FX,Y)\to\Hom_\S((G\comp F)X,GY)
\end{equation} 
for all $X\in\S$ and $Y\in\T$.  Next we observe that for any exact
triangle $X\xto{\a} Y\xto{\b}Z\xto{\g} X[1]$, the map $\b$ is a
retraction if and only if $\g=0$.

The map $F\b_\S$ is not a retraction since $F\g_\S\neq 0$. Thus
$F\b_\S$ factors through $\b_\T$, and $G(F\b_\S)=\b_\S$ factors through
$G\b_\T$. We obtain the following commutative diagram.
$$\xymatrix{ X_{\S} \ar[r]^{\a_{\S}}\ar[d]^{\p} & Y_{\S}
  \ar[r]^{\b_{\S}}\ar[d]^{\psi} & Z_\S \ar[r]^{\g_{\S}}\ar@{=}[d] &
  X_{\S}[1] \ar[d]^{\p[1]}\\ GX_{\T} \ar[r]^{G\a_{\T}} & GY_{\T}
  \ar[r]^{G\b_{\T}} & GZ_\T \ar[r]^{G\g_{\T}}& GX_{\S}[1]}$$
On the
other hand, $G\b_\T$ is not a retraction since the bijection
(\ref{eq:G}) implies $G\g_\T \neq 0$. Thus $G\b_\T$ factors through
$\b_\S$, and we obtain the following commutative diagram.
$$\xymatrix{ GX_{\T} \ar[r]^{G\a_{\T}}\ar[d]^{\p'} & GY_{\T}
\ar[r]^{G\b_{\T}}\ar[d]^{\psi'} & GZ_\T \ar[r]^{G\g_{\T}}\ar@{=}[d] &
GX_{\T}[1] \ar[d]^{\p'[1]}\\ X_{\S} \ar[r]^{\a_{\S}} & Y_{\S}
\ar[r]^{\b_{\S}} & Z_\S \ar[r]^{\g_{\S}}& X_{\S}[1]}$$ We have
$\b_\S\comp (\psi'\comp\psi)=\b_\S$, and this implies that
$\psi'\comp\psi$ is an isomorphism, since $\b_\S$ is right minimal. In
particular, $GY_\T=Y_\S\amalg W$ for some object $W$.  It follows that
$$GX_\T\xto{G\a_\T} GY_\T\xto{G\b_\T}GZ_\T\xto{\g_\T} GX_\T[1]$$ is
the coproduct of $X_\S\xto{\a_\S} Y_\S\xto{\b_\S}Z_\S\xto{\g_\S}
X_\S[1]$ and the triangle $W\xto{\id}W\to 0\to W[1]$.
\end{proof}

Now suppose that $\La$ is an artin algebra. We fix an indecomposable
compact object $Z$ in $\bfK(\Inj\La)$, and we want to compute the
Auslander-Reiten triangle $X\to Y\to Z\to X[1]$. We apply Happel's
functor 
$$H\colon \bfK^c(\Inj\La)\xto{\sim}\bfD^b(\mod\La)\to\umod\hat\La$$
and obtain an indecomposable non-projective $\hat\La$-module
$Z'=HZ$. For instance, if $Z=\bfi N$ is the injective resolution of an
indecomposable $\La$-module $N$, then $H\bfi N=N$ where $N$ is viewed
as a $\hat\La$-module via the canonical algebra homomorphism
$\hat\La\to\La$. Now take the almost split sequence $0\to D\Tr Z'\to
Y'\to Z'\to 0$ in $\Mod\hat\La$. This gives rise to an
Auslander-Reiten triangle $D\Tr Z'\to Y'\to Z'\to D\Tr Z'[1]$ in
$\uMod\hat\La$. We apply the composite
$$\uMod\hat\La\xto{\sim}\bfK_\ac(\Inj\hat\La)
\xto{\Hom_{\hat\La}(\La,-)}\bfK(\Inj\La).$$ It follows from
Proposition~\ref{pr:adj} that the result is a coproduct of the
Auslander-Reiten triangle  $X\to Y\to Z\to X[1]$ and a split exact triangle.

\subsection*{Acknowledgement} 
The authors wish to thank Igor Burban for some helpful discussions on
the topic of this paper. Moreover, we are grateful to Helmut Lenzing
and Dieter Vossieck for pointing out some less well-known references.

\end{document}